\documentclass[11pt]{article}
\usepackage[margin=1in]{geometry}

\usepackage{setspace}
\doublespace

\usepackage{amsmath}
\usepackage{amssymb}
\usepackage[square]{natbib}
\setcitestyle{numbers}

\newtheorem{theorem}{Theorem}[section]
\newtheorem{corollary}{Corollary}[section]
\newtheorem{lemma}{Lemma}[section]
\newtheorem{remark}{Remark}[section]

\newcommand{\qed}{$\square$}                     

\newcommand{\RR}{\mathbb{R}}
\newcommand{\NN}{\mathbb{N}}

\newcommand{\F}{{\mathcal{F}}}

\newcommand{\M}{{\mathcal{M}}}

\newcommand{\ii}{{\mathrm i}}
\newcommand{\ee}{{\mathrm e}}

\newcommand{\dd}{{\mathrm{d}}}

\newcommand{\numberpi}{{\mathrm{\pi}}}

\newcommand{\integrald}{{\mathrm{d}}} %

\newcommand{\dirac}{I}                   

\newcommand{\G}{G}

\newcommand{\Expect}{\mathrm{E}}                 
\newcommand{\Cov}{\mathrm{Cov}}                  

\newcommand{\wE}{\w\Expect}
\newcommand{\wexp}{\w\Expect^{+}}

\newcommand{\exponent}[1]{\exp\{#1\}}           
\newcommand{\Exponent}[1]{\exp\Bigl\{#1\Bigr\}} 
\newcommand{\norm}[1]{\|\, #1\,\|}              
\newcommand{\Norm}[1]{\Big\|\,#1\,\Big\|}       
\newcommand{\ab}[1]{\vert #1\vert}           
\newcommand{\Ab}[1]{\Big\vert #1 \Big\vert}   

\newcommand{\w}{\widehat}

\def\\{\backslash}

\begin{document}


\title{
Smoothing effect of Compound Poisson approximation to distribution
of weighted sums}

\author{V. \v Cekanavi\v cius   and A. Elijio   \\
{\small
Department of Mathematics and Informatics, Vilnius University,}\\
{\small Naugarduko 24, Vilnius 03225, Lithuania.}\\{\small E-mail:
vydas.cekanavicius@mif.vu.lt and aiste.elijio@gmail.com   } }

\date{}

\maketitle

\begin{abstract}

The accuracy of compound Poisson approximation to the sum
$S=w_1S_1+w_2S_2+\dots+w_NS_N$ is estimated.
 Here $S_i$ are sums
of independent or weakly dependent random variables, and $w_i$
denote weights. The overall smoothing effect of $S$ on $w_iS_i$ is
estimated by L\' evy concentration function.
\end{abstract}

\vspace*{.5cm} \noindent {\bf{Key words:}} \emph{\small
characteristic function, concentration function, compound Poisson
distribution, Kolmogorov norm, weighted random variables.}

\vspace*{.5cm} \noindent {\small {\it MSC 2000 Subject
Classification}:
Primary 60F05.   
Secondary 60G50;     
}


\newpage

\section{Introduction}

 Let us consider  typical cluster  sampling design: the
entire population consists of different clusters, and the
probability for each cluster to be selected into the sample is
known.  The sum of sample elements then is equal to
$S=w_1S_1+w_2S_2+\cdots+w_NS_N=w_1(X_{11}+X_{12}+\cdots+X_{1n_1})+\cdots+w_N(X_{N1}+X_{N2}+\cdots
X_{Nn_N})$. Here $w_i$ denote weights, which are inversely
proportional to probabilities to be selected into sample.

We explain  motivating idea of this paper by considering simple
example, when $N=2$ and $w_1=w_2=1$.  We want to estimate
$d(S_1+S_2,Z_1+Z_2)$, where $d(\cdot,\cdot)$ denotes some
probabilistic metric. The majority of metrics allows the following
simplification
\begin{equation}
d(S_1+S_2,Z_1+Z_2)\leqslant d(S_1+S_2,Z_1+S_2)+
d(Z_1+S_2,Z_1+Z_2)\leqslant
d(S_1,Z_1)+d(S_2,Z_2).\label{pirmas}\end{equation}
  Such
approach is reasonable only if both final estimates are of similar
order. Otherwise, by neglecting $S_2$, we can significantly worsen
the overall estimate of the accuracy of approximation. For
example, let $S_1$ have just few summands and  $d(S_1,Z_1)=O(1)$.
Let  $S_2$ have a large number of summands. Then, neither
$S_1+S_2$ nor $Z_1+S_2$ differ much from $S_2$ and, it is  natural
to expect $d(S_1+S_2,Z_1+S_2)$ to be small. If this is the case,
we say that $S_2$ has smoothing effect on $S_1$. Our aim
 is  investigation of such smoothing effects.

Weighting can radically change the structural properties of $S$.
For example, even if all $S_i$ are lattice, the sum $S$ is not
necessarily lattice random variable. Therefore, the standard
approaches (Tsaregradski's inequality, Stein's method) are
inapplicable.

We introduce necessary notation. Let $\F$ (resp.~$\M$) denote the
set of probability distributions (resp.~finite signed measures) on
$\RR$.
  The Dirac measure concentrated at $a$ is denoted by
$\dirac_a$, $\dirac=\dirac_0$. All products and powers of finite
signed measures $W\in\M$ are defined in the convolution sense, and
$W^0=\dirac$. The exponential of $W$ is the finite signed measure
defined by $ \exponent{W}\,=\,\sum_{m=0}^\infty {W^m}/{m!}$. We
denote by $\widehat W(t)$ the Fourier--Stieltjes transform of
$W\in\M$.

The Kolmogorov (uniform) norm $\norm{W}_K$ and the total variation
norm $\norm{W}$ of $W\in\M$ are defined by
\[\norm{W}_K=\sup_{x\in\RR}\ab{W((-\infty,x])},\quad
\norm{W}=W^{+}\{\RR\}+W^{-}\{\RR\},
\]
respectively. Here $W=W^{+}-W^{-}$ is the Jordan-Hahn
decomposition of $W$.  Note that $\norm{W}_K\leqslant\norm{W}$,
$\norm{WV}_K\leqslant \norm{W}\cdot \norm{V}_K$. If $F\in\F$, then
$\norm{F}_K=\norm{F}=1$. For $F\in\F$, $h\geqslant 0$ L\' evy's
concentration function is defined by
\begin{equation*}
Q(F,h)=\sup_xF\{[x,x+h]\}.
\end{equation*}

 All absolute positive
constants are denoted by the same symbol $C$. Sometimes we supply
$C$ with indices. We also assume usual convention $\sum_{j=a}^b=0$
and $\prod_{j=a}^b=1$, if $b<a$.


\section{Known results}

As a rule, the limiting behavior of weighted sums is investigated
with the emphasis on weights,  for example, see \citep{Liang06},
\citep{Rosal98}, \citep{Zhang97}  and the references therein. In
our paper,  emphasis is on the structure of random variables.

Let us assume that all distributions have  finite thee absolute
moments. Then the Berry-Esseen theorem can be used:
\begin{equation}
\Norm{\prod_{i=1}^n F_i-\Phi(\mu, \sigma^2)}_K \leqslant
\frac{C_1\sum_{i=1}^n
\beta_{3i}}{\big(\sum_{i=1}^n\sigma_i^2\big)^{3/2}} . \label{B-E}
\end{equation}
Here $\beta_{3i}$ and $\sigma^2_i$ are the third absolute moment
and variance of $F_i$, respectively. In many cases, the accuracy
in (\ref{B-E}) is of the order $O(n^{-1/2})$. However, this is not
the case when random variables form triangular array and are close
to zero.

 Hipp \citep{Hipp85}  considered smoothing effect in general case of nonnegative random variables
 with some probability mass at zero.
Here we present one improvement of Hipp's result by Roos, which
follows from the more general proposition in \citep{Ro05}. Let all
$B_i$ be concentrated on $(0,\infty)$ and all $p_i<1$, then
\begin{equation}
\Norm{\prod_{i=1}^n
\big((1-p_i)I+p_iB_i\big)-\Exponent{\sum_{i=1}^n p_i(B_i-I)}}_K
\leqslant
\frac{\pi^2}{4}\sum_{i=1}^n\frac{p_i^2}{1-p_i}Q(\widetilde
H,\mu_i). \label{RoosHipp}
\end{equation}
Here $\mu_i=\int x \integrald B_i(x)$ and $\widetilde
H=\Exponent{\sum_{i=1}^n p_i(1-p_i)(B_i-\dirac)/2}$ . The
smoothing effect is estimated by $Q(\widetilde H,\mu_i)$. Note
 that estimate without smoothing effect is equal to $C\min(\sum_1^np_i^2,\max_ip_i)$,
 see Theorem 2.1, p.97 in \citep{AZ88}.

Apart from the  accompanying compound Poisson distribution as in
 (\ref{RoosHipp}) we consider the second
order (signed) compound Poisson approximations, such as
\begin{equation*}
\Exponent{\sum_{i=1}^n\Big(p_i(B_i-\dirac)-p_i^2(B_i-\dirac)^2/2\Big)}.
\label{G2}
\end{equation*}
Analogues of (\ref{RoosHipp}) have been obtained for this
approximation in \citep{Ro05}. For similar approximations see
\citep{BC02},
 \citep{BaXi99},  \citep{Roo03},  and the references
therein.

Note that lower bound estimates of compound Poisson approximation
to weighted sums have been investigated in \citep{CeEl05}.


\section{Results}

\textbf{1. Sums of 1-dependent random variables.} First we
consider the case, when random variables are non-identically
distributed, that is, $S=w_1S_1+w_2S_2+\cdots+w_NS_N$ and
\begin{equation*}
S_m=\sum_{i=1}^{n_m}X_{mi},\quad m=1,\dots,N.
\label{m1}
\end{equation*}
We assume that $S_m$ and $S_j$ are independent when $m\ne j$. On
the other hand, we allow weak dependence of variables in each sum.
Let $X_{m1},X_{m2},\dots,X_{mn_m}$ be 1-dependent.
 We recall that the sequence of
random variables $\{X_j \}_{j \geq 1}$
  is called $k$-dependent if, for $1 < s < t < \infty$, $t- s > m$, the sigma algebras
generated by $X_1,\dots,X_s$ and $X_t,X_{t+1}\dots$ are
independent. Though further on we consider 1-dependent variables,
it is clear that, by grouping consecutive summands, the sum of
$k$-dependent variables can be reduced to the sum of 1-dependent
ones.

We consider the case when all $X_{mk}$ are concentrated at
$0,1,2\dots$. Factorial moments of $X_{mk}$ are defined by
\[\nu_j^{(m)}(k)=\Expect X_{mk}(X_{mk}-1)\cdots(X_{mk}-j+1),\quad
j=1,2,\dots,\quad m=1,\dots,N,\quad k=1,\dots,n_m.
\]
Let
 \[ \Gamma_{m1}=\sum_{k=1}^{n_m}\nu_1^{(m)}(k),\quad
\Gamma_{m2}=\frac{1}{2}\sum_{k=1}^{n_m}[\nu_2^{(m)}(k)-(\nu_1^{(m)}(k))^2]+\sum_{k=2}^{n_m}\Cov(X_{m,k-1},X_{mk}).\]
The distribution of $w_mS_m$ we denote by $F_m$. Next we define
approximating measures:
\begin{eqnarray*}
\Pi_m&=&\exponent{\Gamma_{m1}(\dirac_{w_m}-\dirac)},
\quad\Pi=\prod_{m=1}^N\Pi_m=\Exponent{\sum_{m=1}^N\Gamma_{m1}(\dirac_{w_m}-\dirac)},\label{Pi}\\
\G_m&=&\exponent{\Gamma_{m1}(\dirac_{w_m}-\dirac)+\Gamma_{m2}(\dirac_{w_m}-\dirac)^2},\quad
\G=\prod_{m=1}^N\G_m,\label{G}\\
M_1&=&\Exponent{0.025\sum_{m=1}^N\Gamma_{m1}(\dirac_{w_m}+\dirac_{-w_m}-2\dirac)}.
\end{eqnarray*}
Finally, we define remainder terms. Let   $\wexp(Y_1,Y_2)=\Expect
Y_1Y_2+\Expect Y_1 \Expect Y_2$ and
\begin{eqnarray*}
R_{m0}&=&\sum_{k=1}^n\Big\{\nu_2^{(m)}(k)+(\nu_1^{(m)}(k))^2+\Expect
X_{m,k-1}X_{mk}\Big\},
\label{Rm0}\\
R_{m1}&=&\sum_{k=1}^n\Big\{(\nu_1^{(m)}(k))^3+\nu_1^{(m)}(k)\nu_2^{(m)}(k)+\nu_3^{(m)}(k)\nonumber\\
&&+ [\nu_1^{(m)}(k-2)+\nu_1^{(m)}(k-1)+\nu_1^{(m)}(k)]\Expect
X_{m,k-1}X_{mk}\nonumber\\
&&+\wE^{+}(X_{m,k-1}(X_{m,k-1}-1),X_{mk})+\wE^{+}(X_{m,k-1},X_{mk}(X_{mk}-1))\nonumber\\
&&+\Expect X_{m,k-2}X_{m,k-1}X_{mk}+\Expect X_{m,k-2}\Expect
X_{m,k-1}X_{m,k}\nonumber\\&&+\wE^{+}(X_{m,k-2},X_{m,k-1})\Expect
X_{mk}\Big\}.\label{Rm1}
\end{eqnarray*}

\begin{theorem}\label{1depteorema} Let, for $m=1,2,\dots,N$; $k=1,2,\dots,n_m$,
$\nu_1^{(m)}(k)\leqslant 1/100$, $\nu_2^{(m)}(k)\leqslant\nu_1^{(m)}(k)$, $\quad\nu_3^{(m)}(k)<\infty$ and
\begin{equation}
\sum_{k=1}^{n_m}\nu_2^{(m)}(k)\leqslant\frac{\Gamma_{m1}}{20},\qquad\sum_{k=2}^{n_m}\ab{\Cov(X_{m,k-1},X_{mk})}
\leqslant\frac{\Gamma_{m1}}{20}.\label{3ab}
\end{equation}
 Then, for any
$h>0$,
\begin{eqnarray}
\norm{F-\Pi}_K&\leqslant& C_2 Q(M_1,h)\sum_{m=1}^N R_{m0}\bigg\{
\frac{w_m}{h}\min\bigl(1,\Gamma_{m1}^{-1/2}\bigr)+\min\bigl(1,\Gamma_1^{-1}\bigr)\bigg\},\label{1done}\\
\norm{F-\G}_K&\leqslant& C_3 Q(M_1,h)\sum_{m=1}^NR_{m1}\bigg\{
\frac{w_m}{h}\min\bigl(1,\Gamma_{m1}^{-1}\bigr)+\min\bigl(1,\Gamma_1^{-3/2}\bigr)\bigg\}.\label{1dtwo}
\end{eqnarray}
\end{theorem}

\begin{remark} The choice of approximation in (\ref{1dtwo}) is by
no means restricted to $\G$. For example, let $\Gamma_{m2}>0$.
Then, taking into account Theorem 3.5 and corresponding Lemmas
from \citep{CeVe13} , it is possible reformulate (\ref{1dtwo}) for
the negative binomial approximation.
\end{remark}

As an application  of Theorem \ref{1depteorema} let us consider
weighted sums of 2- runs. Two-runs statistic and its
generalization $k$-runs statistic are one of the best investigated
cases of sums of weakly dependent discrete random variables, see
\citep{BaXi99}, \citep{BrX01}, \citep{Daly12}, \citep{PeCe1},
\citep{WXia08} and the references therein. Let
$X_{mi}=\eta_{mi}\eta_{m,i+1}$, where $\eta_{mi}\sim Be(p_m)$,
$(i=1,2,\dots,n_m+1)$ are independent Bernoulli variables. Then
$S_m$ is the sum of 1-dependent Bernoulli random variables. It is
known that, if $n_m\geqslant 3$, $p_m\leqslant 1/5$, then
\begin{eqnarray}
\Gamma_{m1}&=&np_m^2,\quad\Gamma_{m2}=
\frac{n_mp_m^3(2-3p_m)-2p_m^3(1-p_m)}{2},\quad R_{m1}\leqslant
Cn_mp_m^4,\nonumber\\
 &&\norm{F_m-\G_m}_K\leqslant\norm{F_m-\G_m}\leqslant
C\frac{p_m}{\sqrt{n_m}}, \label{petr}
\end{eqnarray}
 see \citep{PeCe1}. Therefore, the standard application of the triangle inequality as in (\ref{pirmas}) leads to
estimate
\begin{equation}
\norm{F-\G}_K\leqslant C \sum_{m=1}^N\frac{p_m}{\sqrt{n_m}}.
\label{fe}
\end{equation}
 Let us assume that $w_m\asymp C$. If all $p_m$
are sufficiently small, then conditions of Theorem
\ref{1depteorema} are satisfied. Therefore, taking $h=\min w_m/2$
in (\ref{ac4}), we  obtain
\begin{equation}
\norm{F-\G}_K\leqslant C Q(M_1,h)\sum_{m=1}^N
p_m^2\leqslant\frac{C\sum_{m=1}^Np_m^2}{\sqrt{\sum_{m=1}^Nn_mp^2_m}}.\label{oho}
\end{equation}
 Estimate (\ref{oho}) can
be much smaller than  (\ref{fe}). If $p_i=p$, then the smoothing
effect is very obvious:
\begin{equation*}
\norm{F-\G}_K\leqslant
\frac{C(N)p}{\sqrt{n_1+n_2+\cdots+n_N}}\quad \hbox{vs} \quad
\norm{F-\G}_K\leqslant C(N)p
\bigg(\frac{1}{\sqrt{n_1}}+\dots+\frac{1}{\sqrt{n_N}}\bigg).
\end{equation*}
 Note that due to 1-dependence we can not apply (\ref{RoosHipp}).

\textbf{2. Sums of independent random variables satisfying
Franken's condition.} Theorem's \ref{1depteorema} conditions  can
be relaxed if all random variables are independent. Let us
consider  typical case of clustered sample assuming that, in each
sum, all random variables are independent and identically
distributed. More precisely, let, for $m=1,2,\dots,N$, $H_m$ be
concentrated on lattice $0,w_m,2w_m,\dots$, that is,
$H_m=p_{m0}\dirac+p_{m1}\dirac_{w_m}+p_{m2}\dirac_{2w_j}+\dots$.
We denote $j$th factorial moment of $H_m$ by
\[
\nu_j(m)=\sum_{k=0}^\infty k(k-1)\cdots(k-j+1)p_{mk}\] and assume
Franken's condition
\begin{equation}
    \lambda_m:=\nu_1(m)-\nu_1^2(m)-\nu_2(m)>0.
\label{Franken}
\end{equation}
 Franken \citep{Frank64}  proved that, if the main
probabilistic mass of nonnegative integer-valued random variable
is concentrated at zero and unity, then the distribution of sum of
such variables can be approximated by Poisson distribution quite
accurately (see also \citep{Kru86b}). Franken's condition means
that  $\nu_1(m)\leqslant 1$ and $\nu_2(m)\leqslant\nu_1(m)$. It is
much weaker than Theorem's \ref{1depteorema} assumptions
$\nu_1(m)\leqslant 1/100$, $\nu_2(m)\leqslant \nu_1(m)$ and
(\ref{3ab}).

\begin{theorem} \label{treciate} Let $\nu_3(m)<\infty$, $n_m\in\NN$ and let
condition (\ref{Franken}) be satisfied, $(m=1,2,\dots,N).$
 Then, for all $h>0$,
\begin{eqnarray}
\lefteqn{\Norm{\prod_{m=1}^NH_m^{n_m}-\Exponent{\sum_{m=1}^Nn_m\nu_1(m)\big(I_{w_m}-I\big)}}_K
\leqslant C_4 Q(M_2, h)\nonumber}
\hspace{0.2cm}\\
&&\times\sum_{m=1}^N n_m (\nu_{2}(m)+\nu_1^2(m))
\bigg\{\frac{w_m}{h}\min\Big(1,\frac{1}{\sqrt{n_m\lambda_m}}\Big)
+\min\Big(1,\frac{1}{n_m\lambda_m}\Big)\Big(1+\frac{\nu_{1}(m)}{\lambda_m}\Big)\bigg\}
\label{treciate1}
\end{eqnarray}
and
\begin{eqnarray}
\lefteqn{\Norm{\prod_{m=1}^NH_m^{n_m}-\Exponent{\sum_{m=1}^Nn_m\Big(\nu_1(m)\big(\dirac_{w_m}-\dirac\big)+\frac{\nu_2(m)
-\nu_1^2(m)}{2}(\dirac_{w_m}-\dirac)^2\Big)}}_K \nonumber}
\hspace{3cm}\\
&\leqslant& C_5 Q(M_2, h)\sum_{m=1}^N
n_m[\nu_{3}(m)+\nu_1(m)\nu_2(m)+\nu_1^3(m)]\nonumber\\
 &&
\times\bigg\{\frac{w_j}{h}
\min\bigg(1,\frac{1}{n_m\lambda_m}\bigg)
+\min\bigg(1,\frac{1}{(n_m\lambda_m)^{3/2}}\bigg)\bigg(1+\frac{\nu_1(m)}{\lambda_m}\bigg)\bigg\}
. \label{treciate2}
\end{eqnarray}
Here $M_2$ is symmetric  distribution with  $\widehat
M_2(t)=\Exponent{-\sum_{l=1}^N n_l\lambda_l\sin^2(tw_l/2)}$.
\end{theorem}
For any Bernoulli variable Franken's condition is satisfied.
Therefore, assuming $h=\min_m w_m/2$ and applying (\ref{ac4}), we
obtain the following corollary.
\begin{corollary}\label{isvada1} Let
$H_m=(1-p_m)\dirac+p_m\dirac_{w_m}$, $w_m\asymp C$, $p_m\leqslant
C_6<1$, $m=1,\dots,N$. Then
\begin{equation}
\Norm{\prod_{m=1}^NH_m^{n_m}-\Exponent{\sum_{m=1}^Nn_mp_m\big(\dirac_{w_j}-\dirac\big)}}_K
\leqslant C\Big(\sum_{i=1}^N n_ip_i\Big)^{-1/2}
\sum_{m=1}^N\min\{n_mp_m^2,\sqrt{n_m}p_m^{3/2}\}. \label{isvada1a}
\end{equation}
\end{corollary}

It is easy to check, that if $N=n$, $n_j=1$, then
 up
to constant we get  the classical estimate of  Poisson
approximation to the Poisson-binomial distribution with "magic
factor" : $C\sum_1^np_m^2\Big(\sum_1^np_j\Big)^{-1/2}$.

We also can use (\ref{isvada1a}) for comparison to various known
estimates. Let, in Corollary \ref{isvada1}, $N=2$ and
$n_1p_1\geqslant 1$, $n_2p_2\geqslant 1$. Then the estimates in
(\ref{B-E}), (\ref{RoosHipp}) and (\ref{isvada1a}) are of the
order
\[\frac{1}{\sqrt{n_1p_1+n_2p_2}},\quad\frac{n_1p_1^2+n_2p_2^2}{\sqrt{n_1p_1+n_2p_2}},\quad
\frac{p_1\sqrt{n_1p_1}+p_2\sqrt{n_2p_2}}{\sqrt{n_1p_1+n_2p_2}},\]
respectively. Here we used (\ref{ac4}) for upper bound estimate in
(\ref{RoosHipp}). It is easy to check, that the last estimate
always has better order than the second one. Moreover, if $p_1$
and $p_2$ tend to zero sufficiently fast, the last estimate is
sharper than the Berry-Esseen estimate.


\textbf{3. Generalized Poisson-binomial distribution.} We further
relax assumptions on the structure of random variables and
consider the  case when all random variables are independent and
have some probability mass at zero. The supports of random
variables are unnecessary discrete and they might not have any
finite absolute moment apart from the first one. We assume that
random variables in each sum are identically distributed. In
principle, we consider the case similar to the one considered in
 (\ref{RoosHipp}). However, we take an
advantage of the fact that not all distributions are different.
Let $\mu_{m1}=\int_{\RR}\ab{x}B_m\{\dd x\}$ and let $Re\w B_m(t)$
denote the real part of $\w B_m(t)$.

\begin{theorem} \label{antrate} Let  $B_j\in\F$,
$0\leqslant p_j\leqslant \tilde C_7<1$, $\mu_{m1}<\infty$
$(j=1,\dots,N)$. Then, for any $h>0$,
\begin{eqnarray}
\lefteqn{\Norm{\prod_{m=1}^N((1-p_m)\dirac+p_mB_m)^{n_m}-\Exponent{\sum_{m=1}^Nn_mp_m(B_m-\dirac)}}_K
 \nonumber}
\hspace{1.8cm}\\
&\leqslant& C_8 Q(M_3,
h)\sum_{m=1}^Nn_mp_m^2\bigg\{\frac{\mu_{m1}}{h}\min\bigg(1,\frac{1}{\sqrt{n_mp_m}}\bigg)+
\min\bigg(1,\frac{1}{n_mp_m}\bigg)\bigg\}
 \label{antrate1}
\end{eqnarray}
and
\begin{eqnarray}
\lefteqn{\Norm{\prod_{m=1}^N((1-p_m)\dirac+p_mB_m)^{n_m}-\Exponent{\sum_{m=1}^N\big(n_mp_m(B_m-\dirac)-\frac{n_m}{2}p_m^2(B_m-\dirac)^2\big)}}_K
 \nonumber}
\hspace{1cm}\\
&\leqslant& C_9 Q(M_3,
h)\sum_{m=1}^Nn_mp_m^3\bigg\{\frac{\mu_{m1}}{h}\min\bigg(1,\frac{1}{{n_mp_m}}\bigg)+
\min\bigg(1,\frac{1}{(n_mp_m)^{3/2}}\bigg)\bigg\}.
 \label{antrate2}
\end{eqnarray}
Here $M_3$ is symmetric distribution with  $\widehat
M_3(t)=\Exponent{\sum_{l=1}^N0.5n_lp_l(1-p_l)\big(Re\widehat
B_l(t)-1\big)}$.
\end{theorem}

\begin{remark} (i)  Though the accuracy of approximation is
similar to that of previous Theorems, the structure of
approximating Compound Poisson distribution is much more
complicated.

(ii) If,  $B_m\{[0,\infty)\}=1$, $\mu_{m1}\asymp C$,
$(m=1,2,\dots, N)$ then by (\ref{ac4}) we can obtain estimate
similar to (\ref{isvada1a}). Therefore, it is not difficult  to
construct examples similar to the ones, considered for the
previous theorem, and demonstrating the effect of smoothing.

(iii) If $n_j=1$, $N=n$,  then (\ref{antrate1}) is a version of
(\ref{RoosHipp}) for $B_m\{\RR\}=1$. On the other hand, if  all
$B_m\{[0,\infty)\}=1$, then (\ref{RoosHipp}) is more accurate than
(\ref{antrate1}).
\end{remark}


\section{Auxiliary results}

%
Further we  need the following lemmas.
\begin{lemma}\label{ac} Let $F, G \in \F$, $h>0$ and $a>0$. Then
\begin{eqnarray}
Q(F, h)&\leqslant&
\Bigg(\frac{96}{95}\Bigg)^2h\int_{\ab{t}\leqslant 1/h}
\Ab{\widehat F(t)}\,\integrald t, \label{ac1}\\
Q(F, h)&\leqslant& \bigg(1+\Bigg(\frac{h}{a}\Bigg)\bigg) Q(F, a),
\label{ac3}\\
Q(\exponent{a(F-I)}, h)&\leqslant&
\frac{C}{\sqrt{aF\left\{\ab{x}>h\right\} }}. \label{ac4}
\end{eqnarray}
If, in addition, $\widehat F(t)\geqslant 0$, then
\begin{equation}
h \int_{\ab{t}\leqslant 1/h} \ab{\widehat F(t)}\,\integrald t
\leqslant CQ(F, h). \label{ac5}
\end{equation}
\end{lemma}
Lemma \ref{ac} contains  well-known properties of Levy's
concentration function (see, for example, \citep{AZ88}, Chapter
2).

For $h\in(0,\infty)$ and a finite nonnegative measure $G$ on
$\RR$, set $\ab{G}_{h-}=\sup_{x\in\RR}G\{(x,x+h)\}$.
\begin{lemma} \label{vienas} (\citep{CeRo06})
Let $W_1,W_2\in\M$ with $W_1\{\RR\}=0$, and set $W=W_1+W_2$. For
$y\in[0,\infty)$, let
\[\rho(y)=\min\big\{\ab{W^+}_{y-},\,\ab{W^-}_{y-}\big\}.\]
Then, for arbitrary $h\in(0,\infty)$ and $r\in(0,1)$, we have
\[\norm{W}_K\leq
\frac{1}{2r}\,\norm{W_1} +\frac{1}{2\numberpi\,r}\int_{\ab{t}<1/h}
\Ab{\frac{\widehat{W_2}(t)}{t}}\,\integrald t
+\frac{1+r}{2r}\,\rho(4\,\eta(r)h),\] where $\eta(r)\in(0,\infty)$
is defined by the equation
\begin{equation*}\label{eq171}
\frac{1+r}{2}=\frac{2}{\numberpi}\int_0^{\eta(r)}\frac{\sin^2(x)}{x^2}
\,\integrald x.
\end{equation*}
\end{lemma}
\begin{lemma} \label{du} (\citep{CeRo06}) For $F\in\F$, $W\in\M$ with $W\{\RR\}=0$,
and $\vartheta\in(0,\infty)$, we have
\begin{equation}\label{l3}
\ab{(WF)^+}_{\vartheta-}
\leq\frac{1}{2}\,\norm{W}\,\ab{F}_{\vartheta-}.
\end{equation}
\end{lemma}
From Lemmas \ref{vienas} and \ref{du} and (\ref{ac1}) and
(\ref{ac5}) the following result follows
\begin{lemma}\label{ad} Let $h>0$, $W\in \M$, $W\{\RR\}=0$, $P\in \F$, $M$ be distribution with nonnegative characteristic
function and $\ab{\widehat P(t)}\leqslant C M(t)$, for
$\ab{t}\leqslant 1/h$. Then
\begin{eqnarray*}
\norm{W P}_K&\leqslant& C\int_{\ab{t}\leqslant 1/h}
\Ab{\frac{\widehat W(t)\widehat P(t)}{t}}\,\integrald t
+ C\left\|W\right\|Q(P, h)\nonumber\\
&\leqslant& C\Big(\sup_{\ab{t}\leqslant1/h}
\frac{\ab{\widehat W(t)}}{\ab{t}}\cdot\frac{1}{h}+\left\|W\right\|\Big)Q(M,
h).
\label{ad1}
\end{eqnarray*}
\end{lemma}

\textbf{Proof.} We apply Lemma \ref{vienas} with  $W_1=0$,
$W=W_2=WP$, and $r=0.5$. Then by (\ref{ac3}) and (\ref{l3}) we
have
\begin{equation*}
\rho(4\,\eta(r)h)\leqslant
C\ab{(WP)^{+}}_{4\eta(r)h}\leqslant\norm{W}
Q(P,4\eta(r)h)\leqslant C\norm{W}Q(P,h).
\end{equation*}
Moreover, applying (\ref{ac1}) and (\ref{ac5}), we prove that
  \begin{equation*}
Q(P,h)\leqslant Ch\int_{\ab{t}\leqslant 1/h} \ab{\widehat
P(t)}\,\integrald t\leqslant C h\int_{\ab{t}\leqslant 1/h}
\widehat M(t)\,\integrald t\leqslant CQ(M,h)
\end{equation*}
and
  \begin{equation*}
\int_{\ab{t}\leqslant 1/h} \ab{\widehat P(t)}\,\integrald
t\leqslant C\frac{1}{h} h\int_{\ab{t}\leqslant 1/h} \widehat
M(t)\,\integrald t\leqslant \frac{1}{h}CQ(M,h).
\end{equation*}
This, obviously, completes the proof of Lemma. \qed


\begin{lemma}\label{ae} Let $M\in\F$ be concentrated on integers,
$\sum_{k=-\infty}^\infty\ab{kM\{k\}}<\infty$. Then, for all
$\gamma>0$ and $\upsilon\in R$,
\begin{eqnarray*}
\norm{M}^2 \leqslant \Big(\frac{1}{2}+\frac{1}{2\pi\gamma}\Big)
\int_{-\pi}^{\pi} \Big(\gamma\ab{\widehat
M(t)}^2+\frac{1}{\gamma}\Ab{\big(\widehat
M(t)e^{-it\upsilon}\big)^{'}}^2\Big) \,\integrald t. \label{ae1}
\end{eqnarray*}
\end{lemma}
Lemma \ref{ae} has been proved in \citep{Pre85}.

%


\begin{lemma} \label{galbut} Let conditions of Theorem
\ref{1depteorema} be satisfied. Then, for all $t\in\RR$,
$m=1,\dots,N$,
\begin{eqnarray}
\ab{\w F_m(t)}, \ab{\w\G_m(t)},\ab{\w\Pi_m(t)}
&\leqslant&\exponent{-0.26\Gamma_{m1}\sin^2(t_m/2)},\label{glb1}\\
\ab{\w F_m(t)-\w\G_m(t)}&\leqslant&
CR_{m1}\ab{z_m(t_m)}^3\psi_m^{2.6},\label{glb2}\\
\Ab{\Bigl(\exponent{-\ii t_m\Gamma_{m1}}(\w
F_m(t)-\w\G_m(t))\Bigr)'_{t_m}}&\leqslant&CR_{m1}\ab{z(t_m)}^2(1+\ab{z(t_m)}^2\Gamma_{m1})\psi_m^{2.6}
\nonumber\\
&\leqslant& CR_{m1}\ab{z(t_m)}^2\psi_m^2,\label{glb3}\\
\ab{\w
F_m(t)-\w\Pi(t)}&\leqslant&CR_{m0}\ab{z(t_m)}^2\psi_m^{2.6},\\
\Ab{\Bigl(\exponent{-\ii t_m\Gamma_{m1}}(\w
F_m(t)-\w\Pi_m(t))\Bigr)'_{t_m}}&\leqslant&CR_{m0}\ab{z(t_m)}\psi_m^{2}.
\end{eqnarray}
Here $t_m=tw_m$, $z(t_m)=\ee^{\ii t_m}-1$,
$\psi_m=\exponent{-0.1\Gamma_{m1}\sin^2(tw_m/2)}$.
\end{lemma}
 All estimates in Lemma \ref{galbut} follow from  Lemmas 7.4, 7.6, 7.7 and the proof of theorem 5.1 in \citep{CeVe13}.
\section{Proofs} As in previous Section $z(t)=\ee^{\ii
t}-1$, $t_m=tw_m$,
$\psi_m=\exponent{-0.1\Gamma_{m1}\sin^2(tw_m/2)}$. We use the
notation $\theta$ for all quantities satisfying
$\ab{\theta}\leqslant 1$.

\textbf{Proof of Theorem  \ref{1depteorema}.} By properties of the
total variation norm
\begin{eqnarray}
\norm{F-\G}_K&=&\Norm{\prod_{m=1}^NF_m-\prod_{m=1}^N\G_m}_K
\leqslant\sum_{m=1}^N\Norm{(F_m-\G_m)\prod_{l=1}^{m-1}F_l
\prod_{l=m+1}^N\G_l}_K
\nonumber\\
&=&\Norm{(F_m-\G_m)\exponent{-0.05\Gamma_{m1}(\dirac_{w_m}-\dirac)}}_K\nonumber\\
&&\times
\Norm{\exponent{0.05\Gamma_{m1}(\dirac_{w_m}-\dirac)}\prod_{l=1}^{m-1}F_l\prod_{l=m+1}^N\G_l}_K\nonumber\\
&=:&\sum_{m=1}^N\norm{W_mP_m}_K. \label{glb0}
\end{eqnarray}
Note that $\exponent{-0.05\Gamma_{m1}(\dirac_{w_m}-\dirac)}$ is
signed measure of finite variation.

Applying (\ref{glb1}) we obtain
 \begin{equation*}\ab{\w
P_m(t)}\leqslant C\psi_m\prod_{l\ne m}^N\psi_l^{2.6}\leqslant
M_1(t).\label{mHm}\end{equation*}
 Similarly, from (\ref{glb2}) it follows that
 \begin{equation*}
\ab{W(t)}\leqslant
CR_{m1}\ab{z(t_m)}^3\psi_m^{2.6}\psi_m^{-1}\leqslant
CR_{m1}\ab{z(t_m)}^2w_m\ab{t}\psi_m^{1.6}\leqslant
CR_{m1}\min(1,\Gamma_{m1}^{-1})\psi_m^{0.5}w_m\ab{t}. \label{glb4}
 \end{equation*}
Here $\psi_m=\exponent{-0.1\Gamma_{m1}\sin^2(tw_m/2)}$. Applying
Lemma \ref{ad} we obtain
\begin{equation}
\norm{W_mP_m}_K\leqslant
CQ(M_1,h)\bigg\{R_{m1}\min\bigl(1,\Gamma_{m1}^{-1}\bigr)\frac{w_m}{h}+\norm{W_m}\bigg\}.
\label{glb5}
\end{equation}

It remains to estimate $\norm{W_m}$. Since, total variation norm
is invariant to scale change, further we assume $w_m=1$, $t_m=t$.
Then, applying Lemma (\ref{galbut}), we obtain
\begin{eqnarray*}\ab{\w W_m(t)}&\leqslant&
CR_{m1}\min(1.\Gamma_{m1}^{-3/2})\psi^2_m,\\
\Ab{\Big(\exponent{-0.9\ii t\Gamma_{m1}}\w W_m(t)
\Big)'_t}&\leqslant&\Ab{\Bigl(\exponent{-\ii t\Gamma_{m1}}(\w
F_m(t)-\w\G_m(t))\Bigr)'_{t}}\exponent{0.1\Gamma_{m1}\sin^2(t/2)}\\
&&+C\ab{\w
F_m(t)-\w\G_m(t)}\Gamma_{m1}\ab{\sin(t/2)}\exponent{0.1\Gamma_{m1}\sin^2(t/2)}\\
&\leqslant&CR_{m1}\sin^2(t/2)\psi_m^{1.5}(1+\Gamma_{m1}\sin^2(t/2))\leqslant
CR_{m1}\min(1,\Gamma_{m1}^{-1})\psi_m.
\end{eqnarray*}
Taking into account the last two estimates and, applying Lemma
\ref{ae} with $\gamma=\max(1,\sqrt{\Gamma_{m1}})$,
$\upsilon=0.9\Gamma_{m1}$, we get
\[\norm{W_m}\leqslant CR_{m1}\min(1,\Gamma_{m1}^{-3/2}).\]
Substituting the last estimate estimate into (\ref{glb5}) and
(\ref{glb0}) we complete the proof of (\ref{1dtwo}). The proof of
(\ref{1done}) is very similar and, therefore, omitted.  \qed

\textbf{Proof of Theorem \ref{treciate}.} The estimates are proved
exactly by the same arguing as in the proof of Theorem
\ref{1depteorema}. Let
 \begin{eqnarray*}
 \w
D_m(t)&=&\Exponent{\nu_1(m)z(t_m)+ (\nu_2(m)
-\nu_1^2(m))\frac{z^2(t_m)}{2}},\\
\w V_m(t)&=&\bigl(\w H_m^{n_m}(t)-\w
D_m^{n_m}(t)\bigr)\exponent{-0.5n_m\lambda_mz(t_m)},\\
L_m&=& \quad
\Exponent{0.5\lambda_m(\dirac_{w_m}-\dirac)}\prod_{l=1}^{m-1}H_l^{n_l}\prod_{l=m+1}^ND_l^{n_l}.
\end{eqnarray*}
Then
\begin{equation}
\Norm{\prod_{m=1}^NH_m^{n_m}-\prod_{m=1}^ND_m^{n_m}}_K\leqslant\sum_{m=1}^N\norm{V_mL_m}_K.
 \label{Fr1}
\end{equation}
Distribution $H_m$ and approximation $D_m$, for the case $w_m=1$,
were investigated in numerous papers. Let
$r_1(m)=\nu_1^3(m)+\nu_1(m)\nu_2(m)+\nu_3(m)$. Taking into account
Lemmas 2 and 3 in \citep{SC88} and proof of Theorem 3 in
\citep{Kru86b} we can write the following expressions
\begin{eqnarray*}
\lefteqn{\w
H_m(t)=1+\nu_1(m)z(t_m)+\frac{\nu_2(m)}{2}z^2(t_m)+\theta
C\nu_3(m)\ab{z(t_m)}^3,}\hspace{4.5cm}\\
\lefteqn{\big(\w
H_m(t)\big)_{t_m}'=1+\nu_1(m)(z(t_m))'+\frac{\nu_2(m)}{2}(z^2(t_m))'+\theta
C\nu_3(m)\ab{z(t_m)}^2,}\hspace{4.5cm}\\
\ab{\w H_m(t)},\, \ab{\w D_m(t)},
\ab{\exponent{\nu_1(m)z(t_m)}}&\leqslant&\exponent{-2\lambda_m\sin^2(t_m/2)},\\
\ab{\w
H_m(t)-\exponent{\nu_1(m)z(t_m)}}&\leqslant&C(\nu_1^2(m)+\nu_2(m))\ab{z(t_m)}^2,\\
\ab{(\w
H_m(t)-\exponent{\nu_1(m)z(t_m)})_{t_m}'}&\leqslant&C(\nu_1^2(m)+\nu_2(m))\ab{z(t_m)},\\
\ab{\w H_m(t)-\w
D_m(t)}&\leqslant&C r_1(m)\ab{z(t_m)}^3,\\
\ab{(\w H_m(t)-\w D_m(t))_{t_m}'}&\leqslant&C r_1(m)\ab{z(t_m)}^2.
\end{eqnarray*}
%
Therefore, $\ab{L_m(t)}\leqslant C\w M_2(t)$, and
\begin{eqnarray*}
\ab{\w V_m(t)}&\leqslant& Cn_m\ab{\w H_m(t)-\w
D_m(t)}\exponent{n_m\lambda_m\sin^2(t_m/2)-2(n_m-1)\lambda_m\sin^2(t_m/2)}\\
&\leqslant&C\exponent{-n_m\lambda_m\sin^2(t_m/2)}r_1(m)\ab{z(t_m)}^3\\
&\leqslant&C\exponent{-0.5n_m\lambda_m\sin^2(t_m/2)}r_1(m)\min(1,
(n_m\lambda_m)^{-1})w_m\ab{t}.
\end{eqnarray*}
Applying Lemma \ref{ad} we obtain
\begin{equation}
\norm{V_mL_m}_K\leqslant CQ(M_2,h)\bigg\{n_m
r_1(m)\min\bigl(1,(n_m\lambda_m)^{-1}\bigr)\frac{w_m}{h}+\norm{V_m}\bigg\}.
\label{Fr2}
\end{equation}
It remains to estimate $\norm{V_m}$. Total variation norm is
invariant to scale change. Therefore,  we can assume $w_m=1$,
$t_m=t$. For the sake of brevity we use the following notation
omitting the dependence on $t$ and $m$:
\[\omega=\exponent{-0.5n_m\lambda_m\sin^2(t/2)},\quad
u_1=\w H(t)\exponent{-\ii t\nu_1(m)},\quad u_2=\w D(t)\exponent{-
\ii t\nu_1(m)}.\] Taking into account relations from above, we can
write $\ab{\w V_m(t)}\leqslant \omega^2 r_1(m)\min(1,
(n_m\lambda_m)^{-3/2})$ and
\begin{eqnarray*}
\lefteqn{\ab{(\w V_m(t)\exponent{-n_m\nu_1(m)\ii
t+0.5n_m\lambda_m\ii
t})'}=\ab{\big((u_1^{n_m}-u_2^{n_m})\exponent{0.5n_m\lambda_m(1+\ii
t-\ee^{\ii t})}\big)'}}\hspace{2.5cm}\\
&\leqslant&n_m\ab{u_1^{n_m-1}u_1'-u_2^{n_m-1}u_2'}\omega^{-2}+\ab{u_1^{n_m}-u_2^{n_m}}0.5n_m\lambda_m\ab{z(t)}\omega^{-2}\\
&\leqslant&n_m\ab{u_1}^{n_m-1}\ab{u_1'-u_2'}\omega^{-2}+
n_m\ab{u_2'}\ab{u_1^{n_m-1}-u_2^{n_m}}\omega^{-2}\\
&&+\ab{u_1^{n_m}-u_2^{n_m}}0.5n_m\lambda_m\ab{z(t)}\omega^{-2}\\
&\leqslant&n_m\ab{u_1'-u_2'}\omega^2+Cn_m^2(\nu_2+\nu_1)r_1(m)\ab{z(t)}^4\omega^2+Cn_m^2\lambda_m
r_1(m)\ab{z(t)}^4\omega^2\\
&\leqslant&Cn_mr_1(m)\ab{z(t)}^2\omega^{1.5}\bigl(1+n_m[\nu_1(m)+\lambda_m]\ab{z(t)}^2\omega^{0.5}\bigr)\\
&\leqslant&Cn_mr_1(m)\omega\bigg(1+\frac{\nu_1(m)}{\lambda_m}\bigg)\min\bigg(1,\frac{1}{n_m\lambda_m}\bigg).
\end{eqnarray*}
Applying Lemma \ref{ae} with $\gamma=\max(1,\sqrt{n_m\lambda_m})$,
$\upsilon=n_m\nu_1(m)-0.5n_m\lambda_m$, we get
\begin{equation}\norm{V_m}\leqslant
Cn_mr_1(m)\bigg(1+\frac{\nu_1(m)}{\lambda_m}\bigg)\min\bigg(1,\frac{1}{(n_m\lambda_m)^{3/2}}\bigg).\label{Fr4}\end{equation}
Combining the last estimate with (\ref{Fr2}) and (\ref{Fr1}) we
complete the proof of (\ref{treciate2}). The proof of
(\ref{treciate1}) is very similar and , therefore, omitted. \qed

\textbf{Proof of Theorem \ref{antrate}.} Similarly to the proof of
previous Theorem we prove that
\[
\Norm{\prod_{m=1}^N((1-p_m)\dirac+p_mB_m)^{n_m}-
\Exponent{\sum_{m=1}^Np_m(B_m-\dirac)}}_K\leqslant
\sum_{m=1}^N\norm{U_mT_m}_K.
\]
Here \begin{eqnarray*}
 \w U_m(t)&=&\big[((1-p_m)+p_m\w B_m(t))^{n_m}-\exponent{n_mp_m(\w
 B(t)-1)}\big]\\
 &&\times\exponent{0.5n_mp_m(1-p_m)(1-\w B_m(t))},\\
 \w T_m(t)&=&\exponent{0.5n_mp_m(1-p_m)(\w B_m(t)-1)}
 \\
 &&\times\prod_{j=1}^{m-1}((1-p_j)+p_j\w B_j(t))^{n_j}\prod_{j=m+1}^N
 \exponent{n_jp_j(\w B_j(t)-1)}.
\end{eqnarray*}
Taking into account general estimate, $\ab{\w B_m(t)-1}^2\leqslant
2\ab{Re\w B_m(t)-1}$ we obtain:
\begin{eqnarray*}
\ab{\exponent{p_m(\w B_m(t)-1)}}&=&\exponent{p_m(Re\w
B_m(t)-1)},\\
\ab{\exponent{p_m(\w B_m(t)-1)-0.5p_m^2(\w
B_m(t)-1)^2}}&\leqslant&\exponent{p_m(1-p_m)(Re\w B_m(t)-1)},\\
\lefteqn{\ab{1+p_m(\w B_m(t)-1)-\exponent{p_m(\w
B_m(t)-1)}}\leqslant Cp_m^2\ab{\w B_m(t)-1}^2\leqslant
Cp_m^2\ab{Re\w B_m(t)-1}^{1/2}\mu_{m1}\ab{t},}\hspace{9cm}\\
\lefteqn{\ab{1+p_m(\w B_m(t)-1)-\exponent{p_m(\w
B_m(t)-1)-0.5p_m^2(\w B_m(t)-1)^2}}\leqslant Cp_m^3\ab{Re\w
B_m(t)-1}\mu_{m1}\ab{t},}\hspace{8.9cm}.
\end{eqnarray*}
Consequently,
\begin{eqnarray}
\ab{\w
T_m(t)}&\leqslant&\Exponent{0.5\sum_{m=1}^Nn_mp_m(1-p_m)(Re\w
B_m(t)-1)}=\w M_3(t),\label{ros1}\\
\ab{\w U_m(t)}&\leqslant& Cn_m\ab{1+p_m(\w
B_m(t)-1)-\exponent{p_m(\w B_m(t)-1)}}\exponent{0.5p_m(1-p_m)(Re\w
B_m(t)-1)}\nonumber\\
&\leqslant&Cn_mp_m^2\mu_{m1}\ab{t}\min\bigg(1,\frac{1}{\sqrt{n_mp_m}}\bigg).\label{ros2}
\end{eqnarray}
Moreover, due to the properties of total variation norm,
\[\norm{U_m}\leqslant\norm{\bigl((\dirac+p_m(\dirac_1-\dirac))^{n_m}-\exponent{n_mp_m(\dirac_1-\dirac)}
\bigr)\exponent{0.5n_mp_m(1-p_m)(\dirac-\dirac_1)} }.
\]
Arguing similarly as in the proof of (\ref{Fr4}) we prove that
\begin{equation}
\norm{U_m}\leqslant
Cn_mp_m^2\min\bigg(1,\frac{1}{n_mp_m}\bigg).\label{ros3}\end{equation}
From Lemma \ref{ad} and (\ref{ros1})--(\ref{ros3}) we obtain
(\ref{antrate1}). The proof of estimate (\ref{antrate2}) is very
similar and , therefore, omitted.  \qed


\end{document}